\title{ ~~\\ On primes $p$ for which $d$ divides ord$_p(g)$}
\author{Pieter Moree}
\documentclass[12pt]{article}
\usepackage{amssymb, latexsym, amsfonts}
\def\@ptsize{2}
\newtheorem{Thm}{Theorem}

\newtheorem{Lem}{Lemma}

\newtheorem{Cor}{Corollary}

\newtheorem{Prop}{Proposition}
\newcommand{\qed}{\hfill $\Box$}

\begin{document}
\date{}
\maketitle
{\def\thefootnote{}
}

\begin{abstract}
\noindent Let $N_g(d)$ be the set of primes $p$ such that the order of
$g$ modulo $p$, ord$_p(g)$, is divisible
by a prescribed integer $d$. Wiertelak showed
that this set has a natural density, $\delta_g(d)$, with $\delta_g(d)\in \mathbb Q_{>0}$. 
Let $N_g(d)(x)$ be the
number of primes $p\le x$ that are in $N_g(d)$. A simple identity for $N_g(d)(x)$ 
is established. It is used to derive a more compact expression for $\delta_g(d)$ than
known hitherto.
\end{abstract}
\section{Introduction}
Let $g$ be a rational number such that $g\not\in \{-1,0,1\}$ (this assumption on $g$ will
be maintained throughout this note).
Let $N_g(d)$ denote the set of primes $p$ such that the order of $g({\rm mod~}p)$ is divisible
by $d$ (throughout the letter $p$ will also be used to indicate primes). Let $N_g(d)(x)$ denote the 
number of primes in $N_g(d)$ not exceeding $x$. The quantity
$N_g(d)(x)$ (and some variations of it) has been the subject of various publications
[1, 3, 4, 7, 9, 11--19].
Hasse showed that $N_g(d)$ has a Dirichlet density in case $d$ is an odd prime \cite{Hasse1}, respectively
$d=2$ \cite{Hasse2}. The latter case is of
additional interest since $N_g(2)$ is the set of prime divisors of the sequence $\{g^k+1\}_{k=1}^{\infty}$.
(One says that an integer divides a sequence if it divides at least one term of the sequence.)
Wiertelak \cite{W2} established that
$N_g(d)$ has a natural density $\delta_g(d)$ (around the same 
time Odoni \cite{Odoni} did so in the case $d$ is a prime). In a later paper 
Wiertelak \cite{W4} proved, using sophisticated analytic
tools, the following result (with Li$(x)$ the logarithmic integral and with 
$\omega(d)=\sum_{p|d}1$), which gives the best known error term to date.
\begin{Thm} {\rm \cite{W4}}.
We have
$$N_g(d)(x)=\delta_g(d){\rm Li}(x)+O_{d,g}\left({x\over \log^3 x}(\log \log x)^{\omega(d)+1}\right).$$
\end{Thm}
Wiertelak also gave a formula for $\delta_g(d)$ which shows that this is always a positive
rational number.
A simpler formula
for $\delta_g(d)$ (in case $g>0$) has only recently been given by Pappalardi \cite{Pappalardi}.
With some effort Pappalardi's and Wiertelak's expressions can be
shown to be equivalent.\\ 
\indent In this note
a simple identity for $N_g(d)(x)$ will be established (given in Proposition \ref{key}). 
From this it is then inferred that $N_g(d)$ has a natural density $\delta_g(d)$ 
that is given by (\ref{naargraad}), which seems to be the simplest expression involving 
field degrees known for 
$\delta_g(d)$. This expression is then readily evaluated.\\ 
\indent In order to state 
Theorem \ref{main} some
notation is needed. Write $g=\pm g_0^h$, where $g_0$ is
positive and not an
exact power of a rational and $h$ as large as possible. 
Let $D(g_0)$ denote the discriminant of the field $\mathbb Q(\sqrt{g_0})$.
The greatest common divisor of $a$ and $b$ respectively the lowest common multiple
of $a$ and $b$ will be denoted by $(a,b)$, respectively $[a,b]$.
Given an integer $d$, we denote by $d^{\infty}$ the supernatural number (sometimes
called Steinitz number), $\prod_{p|d}p^{\infty}$. Note 
that $(v,d^{\infty})=\prod_{p|d}p^{\nu_p(v)}$.\\

\noindent {\tt Definition}. Let $d$ be even and let $\epsilon_g(d)$ be defined as in Table 1 with
$\gamma={\rm max}\{0,\nu_2(D(g_0)/dh)\}$.

\centerline{{\bf Table 1:} $\epsilon_g(d)$}
\medskip
\begin{center}
\begin{tabular}{|c|c|c|c|}\hline
$g\backslash \gamma$&$\gamma=0$&$\gamma=1$&$\gamma=2$\\
\hline\hline
$g>0$&$-{1/2}$&${1/4}$&${1/16}$\\
\hline
$g<0$&${1/4}$&$-{1/2}$&${1/16}$\\
\hline
\end{tabular}
\end{center}
\medskip
\noindent Note that $\gamma\le 2$. Also note that $\epsilon_g(d)=(-1/2)^{2^{\gamma}}$ if $g>0$.
\begin{Thm}
\label{main}
We have
$$\delta_g(d)={\epsilon_1 \over d(h,d^{\infty})}\prod_{p|d}{p^2\over p^2-1},$$ with
$$\epsilon_1=\cases{1 & if $2\nmid d$;\cr
1+3(1-{\rm sgn}(g))(2^{\nu_2(h)}-1)/4& if $2||d$ and $D(g_0)\nmid 4d$;\cr
1+3(1-{\rm sgn}(g))(2^{\nu_2(h)}-1)/4+\epsilon_g(d)& if $2||d$ and $D(g_0)|4d$;\cr
1 & if $4|d$, $D(g_0)\nmid 4d$;\cr
1+\epsilon_{|g|}(d) & if $4|d$, $D(g_0)|4d$.}
$$
In particular, if $g>0$, then
$$\epsilon_1=\cases{1+(-1/2)^{2^{{\rm max}\{0,\nu_2(D(g_0)/dh)\}}} & if $2|d$ and $D(g_0)|4d$;\cr
1 & otherwise,}$$
and if $h$ is odd, then
$$\epsilon_1=\cases{1+(-1/2)^{2^{{\rm max}\{0,\nu_2(D(g)/dh)\}}} & if $2|d$ and $D(g)|4d$;\cr
1 & otherwise,}$$
\end{Thm}
Using Proposition \ref{key} of Section 2 it is also very easy to infer the following
result, valid under the assumption of the Generalized Riemann Hypothesis (GRH).
\begin{Thm}
\label{mainGRH}
Under GRH we have
$$N_g(d)(x)=\delta_g(d){\rm Li}(x)+O_{d,g}(\sqrt{x}\log^{\omega(d)+1}x),$$
where the implied constant depends at most on $d$ and $g$.
\end{Thm}
\indent In Tables 2 and 3 (Section 6) 
a numerical demonstration of Theorem \ref{main} is given.

\section{The key identity}
Let $\pi_L(x)$ denote the number of unramified primes $p\le x$ that split completely in the
number field $L$. For integers $r|s$ let $K_{s,r}=\mathbb Q(\zeta_s,g^{1/r})$.\\
\indent The starting point of the proof of Theorem \ref{main} is the following proposition. By
$r_p(g)$ the residual index of $g$ modulo $p$ is denoted (we have
$r_g(p)=[\mathbb F_p:\langle g\rangle]$).
Note that ord$_p(g)r_p(g)=p-1$.
\begin{Prop}
\label{key}
We have $N_g(d)(x)=\sum_{v|d^{\infty}}\sum_{\alpha|d}\mu(\alpha)\pi_{K_{dv,\alpha v}}(x)$.
\end{Prop}
{\it Proof}. Let us consider the quantity $\sum_{\alpha|d}\mu(\alpha)\pi_{K_{dv,\alpha v}}(x)$. A prime $p$ counted
by this quantity satisfies $p\le x$, $p\equiv 1({\rm mod~}dv)$ and $r_p(g)=vw$ for some integer
$w$. Write $w=w_1w_2$, with $w_1=(w,d)$. Then the contribution of $p$ to 
$\sum_{\alpha|d}\mu(\alpha)\pi_{K_{dv,\alpha v}}(x)$ is $\sum_{\alpha|w_1}\mu(\alpha)$. We conclude that
$$\sum_{\alpha|d}\mu(\alpha)\pi_{K_{dv,\alpha v}}(x)=$$
\begin{equation}
\label{bais}
\#\{p\le x:p\equiv 1({\rm mod~}dv),~v|r_p(g){\rm ~and~}
({r_p(g)\over v},d)=1\}.
\end{equation}
It suffices to show that
$$N_g(d)(x)=\sum_{v|d^{\infty}}\#\{p\le x:p\equiv 1({\rm mod~}dv),~v|r_p(g){\rm ~and~}
({r_p(g)\over v},d)=1\}.$$
Let $p$ be a prime counted on the right hand side. Note that it is counted only once, namely for
$v=(r_p(g),d^{\infty})$. From ord$_p(g)r_p(g)=p-1$ it is then inferred that $d|{\rm ord}_p(g)$. 
Hence every prime counted on the right hand side is counted on the left hand side as well.
Next consider a prime $p$ counted by $N_g(d)(x)$. It satisfies $p\equiv 1({\rm mod~}d)$. Note there is
a (unique) integer $v$ such that $v|d^{\infty}$, $p\equiv 1({\rm mod~}dv)$ and $(r_p(g)/v,d)=1$. Thus $p$
is also counted on the right hand side. \qed\\

\noindent {\tt Remark} 1. From (\ref{bais}) and Chebotarev's density theorem it follows that
\begin{equation}
\label{klem}
0\le \sum_{\alpha|d}{\mu(\alpha)\over [K_{dv,\alpha v}:\mathbb Q]}\le {1\over [K_{dv,v}:\mathbb Q]}.
\end{equation}

\section{Analytic consequences}
Using Proposition 
\ref{key} it is rather straightforward to establish that $N_g(d)$ has a natural density 
$\delta_g(d)$.
\begin{Lem}
\label{analyse}
Write $g=g_1/g_2$ with $g_1$ and $g_2$ integers. Then
\begin{equation}
\label{pappie}
N_g(d)(x)=\left(\delta_g(d)+O_{d,g}\left({(\log \log x)^{\omega(d)}\over \log^{1/8}x}\right)\right){\rm Li}(x),
\end{equation}
where the implied constant depends at most on $d$ and $g$ and
\begin{equation}
\label{naargraad}
\delta_g(d)=\sum_{v|d^{\infty}}\sum_{\alpha|d}{\mu(\alpha)\over [K_{dv,\alpha v}:\mathbb Q]}.
\end{equation}
\end{Lem}
\begin{Cor}
\label{eentje}
The set $N_g(d)$ has a natural density $\delta_g(d)$.
\end{Cor}
The proof of Lemma \ref{analyse} makes use of the following consequence of the Brun-Titchmarsh
inequality.
\begin{Lem}
\label{brunnie}
Let $\pi(x;l,k)=\sum_{p\le x,~p\equiv l({\rm mod~}k)}1$. Then
$$\sum_{v>z\atop v|d^{\infty}}\pi(x;dv,1)=O_d\left({x\over \log x}{(\log z)^{\omega(d)}\over 
z}\right),$$
uniformly for $3\le z\le \sqrt{x}$.
\end{Lem}
{\it Proof}. On noting that $M_d(x):=\#\{v\le x:v|d^{\infty}\}\le (\log x)^{\omega(d)}/\log 2$, it
straightforwardly follows that
$$\sum_{v>z\atop v|d^{\infty}}{1\over v}=\int_z^{\infty}{dM_d(z)\over z}
\ll_d {(\log z)^{\omega(d)}\over z}.$$ By the Brun-Titchmarsh inequality we have 
$\pi(x;w,1)\ll x/(\varphi(w)\log(x/w))$, where the implied constant is absolute and $w<x$.
Thus 
\begin{equation}
\label{af1}
\sum_{z<v,~dv\le x^{2/3}\atop v|d^{\infty}}\pi(x;dv,1)\ll {x\over \varphi(d)\log x}
\sum_{v>z\atop v|d^{\infty}}{1\over v}\ll_d {x\over \log x}{(\log z)^{\omega(d)}\over z}.
\end{equation} 
Using the trivial estimate $\pi(x;d,1)\le x/d$ we see that
\begin{equation}
\label{af2}
\sum_{dv>x^{2/3}\atop d|v^{\infty}}\pi(x;dv,1)\le \sum_{dv>x^{2/3}\atop v|d^{\infty}}{x\over dv}
\le \sum_{w>x^{2/3}\atop w|d^{\infty}}{x\over w}\ll_d x^{1/3}(\log x)^{\omega(d)}.
\end{equation}
On combining (\ref{af1}) and (\ref{af2}) the proof is readily completed. \qed\\ 

\noindent {\it Proof of Lemma} \ref{analyse}. 
From \cite[Lemma 2.1]{Pappalardi} we recall that there exist absolute constants
$A$ and $B$ such that if $v\le B(\log x)^{1/8}/d$, then
\begin{equation}
\label{pappie2}
\pi_{K_{dv,\alpha v}}(x)={{\rm Li}(x)\over [K_{dv,\alpha v}:\mathbb Q]}+
O_g(x e^{-{A\over dv}\sqrt{\log x}}).
\end{equation}
Let $y=B(\log x)^{1/8}/d$.
From the proof of Proposition \ref{key} we see that
$$N_g(d)(x)=\sum_{v|d^{\infty}\atop v\le y}\sum_{\alpha|d}\mu(\alpha)\pi_{K_{dv,\alpha v}}(x)
+O\left(\sum_{v>y\atop v|d^{\infty}}\pi(x;dv,1)\right)=I_1+O(I_2),$$
say. By Lemma  \ref{brunnie} we obtain that $I_2=O(x(\log\log x)^{\omega(d)}\log^{-9/8}x)$.
Now, by (\ref{pappie2}), we obtain
$$I_1=\sum_{v|d^{\infty}\atop v\le y}\sum_{\alpha|d}{\mu(\alpha)\over [K_{dv,\alpha v}:\mathbb Q]}
+O_{d,g}(y{x\over \log^{5/4}x}).$$
Denote the latter double sum by $I_3$. Keeping in mind Remark 1 we obtain
$$I_3=\delta_g(d)+O\left(\sum_{v|d^{\infty}\atop v> y}\sum_{\alpha|d}{\mu(\alpha)\over 
[K_{dv,\alpha v}:\mathbb Q]}\right).$$
Using (\ref{klem}) and Lemma \ref{degree} it follows that
\begin{eqnarray}
\sum_{v|d^{\infty}\atop v> y}\sum_{\alpha|d}{\mu(\alpha)\over [K_{dv,\alpha v}:\mathbb Q]}
&=& O\left(\sum_{v|d^{\infty}\atop v>y}{1\over [K_{dv,v}:\mathbb Q]}\right)
=O({1\over \varphi(d)}\sum_{v|d^{\infty}\atop v>y}{h\over v^2})\nonumber\\
&=&O_d({h(\log y)^{\omega(d)}\over y})=O_{d,g}\left({(\log y)^{\omega(d)}\over y}\right),\nonumber
\end{eqnarray}
and hence $$I_3=\delta_g(d)+O_{d,g}\left({(\log y)^{\omega(d)}\over y}\right).$$
The result follows on collecting the various estimates. \qed
\section{The evaluation of the density $\delta_g(d)$}
A crucial ingredient in the evaluation of $\delta_g(d)$ is the following lemma. 
\begin{Lem} {\rm \cite{Moreealleen}}.
\label{degree}
Write $g=\pm g_0^h$, where $g_0$ is
positive and not an
exact power of a rational. 
Let $D(g_0)$ denote the discriminant of the field $\mathbb Q(\sqrt{g_0})$.
Put $m=D(g_0)/2$ if $\nu_2(h)=0$ and $D(g_0)\equiv 4({\rm mod~}8)$
or $\nu_2(h)=1$ and $D(g_0)\equiv 0({\rm mod~}8)$, and
$m=[2^{\nu_2(h)+2},D(g_0)]$ otherwise. 
Put $$n_r=\cases{m &if $g<0$ and $r$ is odd;\cr
[2^{\nu_2(hr)+1},D(g_0)] &otherwise.}$$
We have
$$[K_{kr,k}:\mathbb Q]=[\mathbb Q(\zeta_{kr},g^{1/k}):\mathbb Q]={\varphi(kr)k\over
\epsilon(kr,k)(k,h)},$$
where, for $g>0$ or $g<0$ and $r$ even we have
$$\epsilon(kr,k)=\cases{2 &if $n_r|kr$;\cr
1 &if $n_r\nmid kr$,}$$
and for $g<0$ and $r$ odd we have
$$\epsilon(kr,k)=\cases{2 &if $n_r|kr$;\cr
 {1\over 2} &if $2|k$ and $2^{\nu_2(h)+1}\nmid k$;\cr
1 &otherwise.}$$
\end{Lem}
{\tt Remark}. Note that if $h$ is odd, then $n_r=[2^{\nu_2(r)+1},D(g)]$. Note that
$n_r=n_{\nu_2(r)}$.\\

\noindent The `generic' degree of $[K_{dv,\alpha v}:\mathbb Q]$ equals 
$\varphi(dv)\alpha v/(\alpha v,h)$ and on
substituting this value in (\ref{naargraad}) we obtain the quantity $S_1$ which is
evaluated in the following lemma.
\begin{Lem}
\label{basis}
We have
$$S_1:=\sum_{v|d^{\infty}}\sum_{\alpha|d}{\mu(\alpha)(\alpha v,h)\over 
\varphi(dv)\alpha v}
=S(d,h),$$
where $$S(d,h):={1\over d(h,d^{\infty})}\prod_{p|d}{p^2\over p^2-1}.$$
\end{Lem}
{\it Proof}. Since for $v|d^{\infty}$ we have $\varphi(dv)=v\varphi(d)$, we can write
$$S_1={1\over \varphi(d)}\sum_{v|d^{\infty}}\sum_{\alpha|d}{\mu(\alpha)(\alpha v,h)\over \alpha v^2}
={1\over \varphi(d)}\sum_{v|d^{\infty}}{(v,h)\over v^2}
\sum_{\alpha|d}{\mu(\alpha)(\alpha v,h)\over \alpha(v,h)}.$$
The expression in the inner sum is multiplicative in $\alpha$ and hence
$$\sum_{\alpha|d}{\mu(\alpha)(\alpha v, h)\over \alpha(v,h)}
=\prod_{p|d}\left(1-{(pv,h)\over p(v,h)}\right)=\cases{{\varphi(d)\over d} & if $(h,d^{\infty})|(v,d^{\infty})$; \cr
0 & otherwise.}$$
On noting that $(v,h)/v^2$ is multiplicative in $v$ and that for $k\ge \nu_p(h)$
$$\sum_{r=k}^{\infty}{(p^r,h)\over p^{2r}}={p^{\nu_p(h)+2-2k}\over p^2-1},$$
one concludes that
$$S_1={1\over d}\sum_{v|d^{\infty}\atop 
(h,d^{\infty})|v}{(v,h)\over v^2}={1\over d}
\prod_{p|d}\sum_{r\ge \nu_p(h)}{(p^r,h)\over p^{2r}}={1\over d}\prod_{p|d}{p^{2-\nu_p(h)}\over p^2-1}=S(d,h).$$
This completes the proof. \qed\\

\noindent {\tt Remark}. Note that the condition $(h,d^{\infty})|(v,d^{\infty})$ is equivalent 
with $\nu_p(v)\ge \nu_p(h)$ for
all primes $p$ dividing $d$.\\

\noindent By a minor modification of the proof of the latter result we infer:
\begin{Lem}
\label{stweek}
Let $k\ge 0$ be an integer. Then
$$S_2(k):=\sum_{v|d^{\infty}\atop \nu_2(v)\ge 
\nu_2(h)+k}\sum_{\alpha|d}{\mu(\alpha)(\alpha v,h)\over \varphi(dv)\alpha v}
=4^{-k}S(d,h).$$
\end{Lem}
The next lemma gives an evaluation of yet another variant of $S_1$.
\begin{Lem}
\label{sdrie}
Let $D$ be a fundamental discrimant. 
Then
$$S_3(D):=
\sum_{v|d^{\infty}\atop 
[2^{\nu_2(hd/\alpha)+1},D]|dv}\sum_{\alpha|d}{\mu(\alpha)(\alpha v,h)\over \varphi(dv)\alpha v}
=\cases{4^{-\gamma}S(d,h) & if $2|d,~D|4d$ and $\gamma\ge 1$;\cr
-{S(d,h)\over 2} & if $2|d,~D|4d$ and $\gamma=0$; \cr
0 & otherwise,}$$
where $\gamma={\rm max}\{0,\nu_2(D/dh)\}$.
\end{Lem}
{\it Proof}. The integer $[2^{\nu_2(hd/\alpha)+1},D]$ is even and is required to divide $d^{\infty}$, hence
$S_3(D)=0$ if $d$ is odd. Assume that $d$ is even. 
If $D$ has an odd prime divisor not divinding $d$, then $D\nmid d^{\infty}$ and hence $S_3(D)=0$.
On noting that $\nu_2(D)\le \nu_2(4d)$ and that the odd part of $D$ is 
squarefree, it follows that if $S_3(D)\ne 0$, then $D|4d$. So assume that $2|d$ and $D|4d$. Note 
that the condition
 $[2^{\nu_2(hd/\alpha)+1},D]|dv$ is equivalent with 
 $\nu_2(v)\ge \nu_2(h)+{\rm max}\{1,\nu_2(D/dh)\}$  for the $\alpha$ that are odd, and
 $\nu_2(v)\ge \nu_2(h)+\gamma$ for the even $\alpha$. Thus if $\gamma\ge 1$ the condition 
 $[2^{\nu_2(hd/\alpha)+1},D]|dv$ is equivalent with $\nu_2(v)\ge \nu_2(h)+\gamma$ and
 then, by Lemma \ref{stweek}, $S_3(D)=S_2(\gamma)=4^{-\gamma}S(d,h)$. If $\gamma=0$ then
 $$S_3(D)=S_2(0)-\sum_{v|d^{\infty}\atop 
\nu_2(v)=\nu_2(h)}\sum_{\alpha|d\atop 2\nmid \alpha}{\mu(\alpha)(\alpha v,h)\over \varphi(dv)\alpha v}.$$
By Lemma \ref{stweek} it follows that $S_2(0)=S(d,h)$. A variation of Lemma \ref{basis} yields that
the latter double sum equals $3S(d,h)/2$. \qed\\

\noindent {\tt Remark}. Put
$$\epsilon_2(D)=\cases{(-1/2)^{2^{{\rm max}\{0,\nu_2(D/dh)\}}} &if $2|d$ and $D|4d$;\cr
0 &otherwise.}$$ 
Note that Lemma \ref{sdrie} can be rephrased as stating that if $D$ is a fundamental discriminant,
then $S_3(D)=\epsilon_2(D) S(d,h)$.\\

\noindent Let $g>0$.
It turns out that ord$_p(g)$ is very closely related to ord$_p(-g)$ and this
can be used to express $N_{-g}(d)(x)$ in terms of $N_g(*)(x)$. From
this $\delta_{-g}(d)$ is then easily evaluated, once one has evaluated $\delta_g(d)$.
\begin{Lem}
\label{zeven}
Let $g>0$. Then
$$N_{-g}(d)(x)=\cases{N_g({d\over 2})(x)+N_g(2d)(x)-N_g(d)(x)+O(1) & if $d\equiv 2({\rm mod~}4)$;\cr
N_g(d)(x)+O(1)& otherwise.}$$
In particular, 
$$\delta_{-g}(d)=\cases{\delta_g({d\over 2})+\delta_g(2d)-\delta_g(d) &if $d\equiv 2({\rm mod~}4)$;\cr
\delta_g(d) & otherwise.}$$
\end{Lem}
The proof of this lemma is a consequence of Corollary \ref{eentje} and the following observation.
\begin{Lem}
Let $p$ be odd and $g\ne 0$ be a rational number. Suppose that $\nu_p(g)=0$. Then
$${\rm ord}_p(-g)=\cases{2{\rm ord}_p(g) & if $2\nmid {\rm ord}_p(g)$;\cr
{\rm ord}_p(g)/2 & if ${\rm ord}_p(g)\equiv 2({\rm mod~}4)$;\cr
{\rm ord}_p(g) & if $4|{\rm ord}_p(g)$.}$$
\end{Lem}
{\it Proof}. Left to the reader. \qed\\

\noindent {\tt Remark}. It is of course also possible to evaluate $\delta_g(d)$ for negative $g$ using
the expression (\ref{naargraad}) and Lemma \ref{degree}, however, this turns out to be
rather more cumbersome than proceeding as above.

\section{The proofs of Theorems \ref{main} and \ref{mainGRH}}
\noindent {\it Proof of Theorem} \ref{main}. 
By Lemma \ref{analyse} it suffices to show that
$$\sum_{v|d^{\infty}}\sum_{\alpha|d}{\mu(\alpha)\over [K_{dv,\alpha v}:\mathbb Q]}=\epsilon_1 S(d,h)$$
If $g>0$, then it follows by Lemma
\ref{degree} that $\delta_g(d)=S_1+S_3(D(g_0))$ and by Lemmas \ref{basis} and \ref{sdrie}
(with $D=D(g_0)$), the claimed evaluation then results in this case. If $h$ is odd, then similarly,
$\delta_g(d)=S_1+S_3(D(g))$ (cf. the remark following Lemma \ref{degree}) 
and, again by Lemma \ref{basis} and \ref{sdrie}, the claimed evaluation then is deduced in this case.
If $g<0$, the result follows after some computation on invoking Lemma \ref{zeven} 
and the result for $g>0$. \qed\\

\noindent {\it Proof of Theorem} \ref{mainGRH}. Recall that $\pi_L(x)$ denotes the number of unramified primes $p\le x$ that split completely in the
number field $L$. Under GRH it is known, cf. \cite{L}, that
$$\pi_L(x)={{\rm Li}(x)\over [L:\mathbb Q]}+
O\left({\sqrt{x}\over [L:\mathbb Q]}\log(d_Lx^{[L:\mathbb Q]})\right),$$
where $d_L$ denotes the absolute discriminant of $L$. From this it follows 
on using the estimate 
$\log |d_{K_{dv_1,\alpha v}}|\le dv(\log (dv)+\log |g_1 g_2|)$ from \cite{Moreealleen} that, 
uniformly in $v$,
$$\pi_{K_{dv,\alpha v}}(x)={{\rm Li}(x)\over [K_{dv,\alpha v}:\mathbb Q]}+O_{d,g}(\sqrt{x}\log x),$$
where $\alpha$ is an arbitrary divisor of $d$.
On noting that in Proposition \ref{key} we can restrict to those integers $v$ satisfying
$dv\le x$ and hence the number of non-zero terms in Proposition \ref{key} is bounded
above by $2^{\omega(d)}(\log x)^{\omega(d)}$, the result easily follows. \qed

\section{Some examples}
In this section we provide some numerical demonstration of our results.\\
\indent The numbers in the column `experimental' arose on counting how many
primes $p\le 
p_{10^8}=2038074743$ with $\nu_p(g)=0$, 
satisfy $d|{\rm ord}_p(g)$.\\

\centerline{{\bf Table 2:} The case $g>0$}
\begin{center}
\begin{tabular}{|c|c|c|c|c|c|c|c|c|}\hline
$g$&$g_0$&$h$&$D(g_0)$&$d$&$\epsilon_1$&$\delta_g(d)$&numerical&experimental\\
\hline\hline
2&2&1&8&2&$17/16$&$17/24$&$0.70833333\cdots$&0.70831919\\ \hline
2&2&1&8&4&$5/4$&$5/12$&$0.41666666\cdots$&0.41667021\\ \hline
2&2&1&8&8&$1/2$&$1/12$&$0.08333333\cdots$&0.08333144\\ \hline
3&3&1&12&11&1&$11/120$&$0.09166666\cdots$&0.09165950\\ \hline
3&3&1&12&12&$1/2$&$1/16$&$0.06250000\cdots$&0.06249098\\ \hline
4&2&2&8&5&1&$5/24$&$0.20833333\cdots$&0.20833328\\ \hline
4&2&2&8&6&$5/4$&$5/32$&$0.15625000\cdots$&0.15625824\\ \hline
\end{tabular}
\end{center}

\centerline{{\bf Table 3:} The case $g<0$}
\medskip
\begin{center}
\begin{tabular}{|c|c|c|c|c|c|c|c|c|}\hline
$g$&$g_0$&$h$&$D(g_0)$&$d$&$\epsilon_1$&$\delta_g(d)$&numerical&experimental\\
\hline\hline
-2&3&1&8&2&$17/16$&$17/24$&$0.70833333\cdots$&0.70835101\\ \hline
-2&2&1&8&4&$5/4$&$5/12$&$0.41666666\cdots$&0.41667021\\ \hline
-2&2&1&8&6&$17/16$&$17/64$&$0.26562500\cdots$&0.26562628\\ \hline
-3&3&1&12&5&1&$5/24$&$0.20833333\cdots$&0.20834107\\ \hline
-3&3&1&12&12&$1/2$&$1/16$&$0.06250000\cdots$&0.06249098\\ \hline
-4&2&2&8&2&2&$2/3$&$0.66666666\cdots$&0.66666122\\ \hline
-4&2&2&8&4&$1/2$&$1/8$&$0.08333333\cdots$&0.08333144\\ \hline
-9&3&2&12&2&$5/2$&$5/6$&$0.83333333\cdots$&0.83333215\\ \hline
-9&3&2&12&6&$11/4$&$11/32$&$0.34375000\cdots$&0.34375638\\ \hline

\end{tabular}
\end{center}

\medskip

\noindent {\bf Acknowledgement}. I like to thank Francesco Pappalardi for
sending me his paper \cite{Pappalardi}. Theorem 1.3 in that
paper made me realize that a relatively simple formula for
$\delta_g(d)$ exists. The data in the tables are produced by a $C^{++}$ program
kindly written by Yves Gallot.

{\small

\medskip\noindent {\footnotesize Max-Planck-Institut f\"ur Mathematik,
Vivatsgasse 7, D-53111 Bonn, Germany.\\
e-mail: {\tt moree@mpim-bonn.mpg.de}}

\end{document}